\documentclass[12pt]{article}
\usepackage{amsfonts}
\usepackage{amssymb}
\usepackage{amsmath}
\usepackage{eufrak}
\usepackage[dvips]{graphicx}
\usepackage{latexsym}
\newcommand{\ind}{\makebox[1em]{\raisebox{-.5ex}[0ex][0ex]{\makebox[0em]%
{$\smile$}}\raisebox{.4ex}[0ex][0ex]{\makebox[-.02em]{$|$}}}}

\newcommand{\nmind}{\makebox[1em]{\raisebox{1.5ex}[0ex][0ex]{\makebox[0em]%
{$\scriptscriptstyle n\! m$}}\makebox[-1em]{$\ind$}}}

\newcommand{\C}{{\EuFrak C}}

\newcommand{\bi}{\begin{itemize}}
\newcommand{\ei}{\end{itemize}}

\newtheorem{theorem}{Theorem}[section]
\newtheorem{lemma}[theorem]{Lemma}
\newtheorem{fact}[theorem]{Fact}

\newtheorem{corollary}[theorem]{Corollary}

\newtheorem{remark}[theorem]{Remark}
\newtheorem{conjecture}[theorem]{Conjecture}
\newtheorem{question}[theorem]{Question}

\title{On relationships between algebraic properties of groups and rings in some model-theoretic contexts}
\author{Krzysztof Krupi\'nski\footnote{Research supported by the Polish Government grant N N201 545938}}
\date{}

\begin{document}
\maketitle
\begin{abstract}
We study relationships between certain algebraic properties of groups and rings definable in a first order structure or $*$-closed in a compact $G$-space. 
%in the context of $\omega$-categorical structures as well as in the context of small Polish structures. For example, we prove that each solvable, $\omega$-categorical group such that every ring interpretable in it is nilpotent-by-finite must be nilpotent-by-finite.
As a consequence, we obtain a few structural results about $\omega$-categorical rings as well as about small, $nm$-stable compact $G$-rings, and we also obtain surprising relationships between some conjectures concerning small profinite groups.  
\end{abstract}

\section[\mbox{}]{Introduction}

One of the main goals in model theory is to classify interesting algebraic structures satisfying some natural model-theoretic assumptions. The algebraic structures that we concentrate on are groups and rings. Our model-theoretic/topological contexts will be mainly $\omega$-categorical structures and small Polish structures (in particular, small profinite structures). In all these contexts, there are many results and conjectures describing the structure of groups or rings satisfying some extra model-theoretic assumptions. However, it seems that relationships between results and conjectures about groups and rings have not been fully described. So, the main goal of this paper is to analyze such relationships.
 
%We will prove that some statements about rings are equivalent to the corresponding statements about groups. 
In Section 2, we study such relationships in a very general context of groups and rings definable in an arbitrary first order structure or $*$-closed in an arbitrary compact $G$-space. In both these contexts, we transfer some properties of groups (e.g. virtual commutativity or solvability) to the corresponding properties of rings (e.g. virtual triviality or nilpotency). As a consequence, using certain results on groups, in Section 3 we get structural results on $\omega$-categorical rings and on small, $nm$-stable compact $G$-rings. For example, we prove that each small, $nm$-stable compact $G$-ring is nilpotent-by-finite. We also obtain surprising relationships between some conjectures on small profinite groups. In our investigations, besides results  proved in Section 2, we will also need to use \cite[Theorem 3.5]{KW} which allows one to deduce certain property of groups (virtual nilpotency of solvable small profinite groups) from the corresponding property of rings (virtual nilpotency of small profinite rings).  Theorem 3.5 from \cite{KW} is restricted to the context of small profinite structures. At the end of the paper, we prove a counterpart of this result for $\omega$-categorical structures.

\section{Preliminaries}

Recall that a first order structure $M$ in a countable language is said to be $\omega$-categorical if, up to isomorphism, $Th(M)$ has at most one model of cardinality $\aleph_0$.

Now, we recall some notions concerning Polish structures. For more details see \cite{Kr}. 
A Polish structure is a pair $(X,G)$, where $G$ is a Polish group acting faithfully on a set $X$ so that the stabilizers of all singletons are closed subgroups of $G$. We say that $(X,G)$ is small if for every $n \in \omega$, there are only countably many orbits on $X^n$ under the action of $G$. 
In a small Polish structure $(X,G)$, there is a ternary relation, $\nmind$, on finite tuples (or subsets) of $X$ which has some nice properties (e.g. symmetry, transitivity, the existence of independent extensions), and which allows us to define the so-called ${\cal NM}$-rank on orbits of finite tuples over finite subsets of $X$ in analogy with the Lascar U-rank. A Polish structure is said to be $nm$-stable if ${\cal NM}$-rank of every orbit is an ordinal.  For $Y\subseteq X^n$, we define $Stab(Y):=\{ g \in G: g[Y]=Y\}$.  We say that $Y$ is invariant [over a finite set $A$] if $Stab(Y)=G$ [$Stab(Y) \supseteq G_A$ where $G_A$ is the pointwise stabilizer of $A$, respectively].

For us, compact spaces and topological groups are Hausdorff by definition.
A [compact] $G$-space is a Polish structure $(X,G)$, where $G$ acts continuously on a [compact] space $X$.
If $(X,G)$ is a $G$-space, we say that $D\subseteq X^n$ is $A$-closed (for a finite $A\subseteq X$) if it is closed and invariant over $A$. We say that it is $*$-closed, if it is $A$-closed for some $A$. Assume $(X,G)$ is a compact $G$-space. We define $X^{teq}$ (topological imaginary extension) as the disjoint union of the spaces $X^n/E$ with $E$ ranging over all $\emptyset$-closed equivalence relations on $X^n$. The spaces $X^n/E$ are called topological sorts of $X^{teq}$. Then, each topological sort $X/E$ together with the group $G/G_{X/E}$ is a compact $G/G_{X/E}$-space. If $E$ is $A$-closed for some finite set $A$, then replacing $G$ by $G_A$, $X/E$ can also be treated as a topological sort.

Now, we recall some facts about groups. A [compact] $G$-group [$G$-ring]
is a Polish structure $(H,G)$, where $G$ acts continuously and by automorphisms on a [compact] topological group [ring] $H$. 
We should mention here that by \cite[Proposition 5.1.2]{RZ}, for a topological ring $R$, the following conditions are equivalent: (i) $R$ is compact, (ii) $R$ is profinite, (iii) there is a basis of open neighborhoods of $0$ consisting of open ideals.  

Let $(X,G)$ be a $G$-space. We say that a group $H$ is a $*$-closed in $X$ [or in $X^{teq}$ if $X$ is compact] if both $H$ and the group operation on $H$ are $*$-closed in $X$ [$X^{teq}$, respectively]. $*$-closed rings are defined analogously.

We work in a compact $G$-space $(X,G)$.
Assume $H$ is $*$-closed in $X^{teq}$; for simplicity, $\emptyset$-closed (then $(H,G/G_H)$ is a compact $G/G_H$-group). Let $a \in H$ and $A\subseteq X$ be finite. We say that the orbit $o(a/A)$ is $nm$-generic (or that $a$ is $nm$-generic over $A$) if for all $b \in H$ with $a\nmind_A b$, one has that $b\cdot a\nmind A,b$. It turns out that under the assumption of smallness of $(X,G)$, $nm$-generics satisfy all the basic properties that generics satisfy in simple groups, including %%%Frank
existence. More precisely, \cite[Proposition 5.5]{Kr} tells us that $o(a/A)$ is $nm$-generic iff $o(a/A) \subseteq_{nm} H$, and by smallness such an orbit exists. 
%Using this, we get the following Lascar inequalities for groups and corollary (see \cite[Corollary 5.6]{kr}). Recall that the $\NM$-rank of an $A$-closed subset $D$ of $X^{teq}$ is the supremum of $\NM(d/A)$, $d \in D$. 

Profinite spaces will be always inverse limits of countable systems.
A profinite structure is a compact $G$-space $(X,G)$, where $X$ is a  profinite space and $G$ is a compact (equivalently profinite) group.
A profinite group [ring] regarded as profinite structure is a compact $G$-group [$G$-ring] $(H,G)$, where $H$ is a profinite group [ring] and $G$ is a profinite group.
Originally profinite structures were defined as pairs $(X,G)$, where $H$ is a profinite space and $G$ is a closed subgroup of the group of all homeomorphisms of $X$ preserving a distinguished inverse system defining $X$; similarly for profinite groups and rings.  By \cite[Proposition 1.5]{Kr2} together with \cite[Proposition 1.4]{Ne} and \cite[Remark 2.9]{KW}, both versions of the  definitions are equivalent.
This means that a profinite group [ring] regarded as profinite structure has a basis of open neighborhoods of $e$ [$0$] consisting of clopen normal subgroups [ideals] invariant under $G$.

It was noticed in \cite{Kr} that if $(H,G)$ is a small compact $G$-group, then $H$ is locally finite and hence profinite. However, $G$ is only Polish (not necessarily compact), which makes the class of small compact $G$-groups much wider than the class of small profinite groups.

Recall some basic notions from ring theory. In this paper, all rings are associative, but they are not assumed to contain 1  or to be commutative. An element $r$ of a ring $R$ is nilpotent of nilexponent $n$ if $r^n=0$, and $n$ is the smallest number with this
property. The ring is nil [of nilexponent $n$] if every element is nilpotent [of nilexponent $\leq n$, and there is an element of nilexponent $n$]. The  ring is nilpotent of class $n$ if $r_1\cdots r_n=0$ for all $r_1,\ldots,r_n\in R$, and $n$ is the smallest number with this property. An element $r$ is null if $rR=Rr=\{0\}$. The ring is null if all its elements are. If $S \subseteq R$, then $Ann_R(S) =\{ r \in R: rS=Sr=\{ 0 \} \}$ is the two-sided annihilator of $S$ in $R$. Note that $Ann_R(S)$ is always a subgroup of $R^+$, and if $S$ is an ideal of $R$, then so is $Ann_R(S)$. $Ann_R^{left}(S)$ and $Ann_R^{right}(S)$ will denote the left and the right annihilator of $S$ in $R$. By $R^n$ we will denote the subring of $R$ generated by all products $r_1\dots r_n$ for $r_1,\dots,r_n \in R$. Notice that $R^n$ is always an ideal.

For groups [rings] $G$ and $H$, $H<G$ means that $H$ is a (not necessarily proper) subgroup [subring] of $G$.

Recall that an abstract group is said to be abelian-by-finite [nilpotent-by-finite, solvable-by-finite] if it has a (normal) abelian [nilpotent, solvable] subgroup of finite index. Similarly, a ring is null-by-finite [nilpotent-by-finite] if it has a null [nilpotent] ideal of finite index. \cite[Lemma 1]{Le} says that if $S$ is a finite index subring of a ring $R$, then $R$ has a finite index ideal contained in $S$. Thus, we have

\begin{remark}\label{Lewin}
A ring is null-by-finite [nilpotent-by-finite] iff it has a null [nilpotent] subring of finite index. 
\end{remark}

We say that a group $G$ is (finite central)-by-abelian-by-finite if there are $H \lhd G$ and $K \lhd H$ such that $G/H$ is finite, $H/K$ is abelian, and $K$ is a finite subgroup of $Z(H)$. A ring $R$ is (finite null)-by-null-by-finite if there are $I\lhd R$ and $J \lhd I$ such that $R/I$ is finite, $I/J$ is null, and $J$ is a finite subring of $Ann_I(I)$.

There are many results describing the structure of $\omega$-categorical groups and small compact $G$-groups. Here, we recall only a few of them, which will be useful in this paper. The first one is \cite[Theorem 1.2]{Ma}, and the second one is \cite[Theorem 5.19]{Kr}.

\begin{fact}\label{Mac NSOP}
An $\omega$-categorical group satisfying NSOP (the negation of the strict order property) is nilpotent-by-finite.
\end{fact}

\begin{fact}\label{nm-stable are solvable}
A small, $nm$-stable compact $G$-group is solvable-by-finite.
\end{fact}

In fact, using the above result, it was shown in \cite{KW2} that small, $nm$-stable compact $G$-groups are nilpotent-by-finite.
It is conjectured that they are even abelian-by-finite. A partial result in this direction is \cite[Theorem 2.9]{KW2}

\begin{fact}\label{NM finite implies abelian}
Let $(H,G)$ be a small compact $G$-group. If ${\cal NM}(H)<\omega$ or ${\cal NM}(H)=\omega^{\alpha}$ for some ordinal $\alpha$, then $H$ is abelian-by-finite.
\end{fact}

\section{General context}

In this section, ${\mathcal M}$ is a first order structure. First, we show that if all groups definable in ${\mathcal M}$ have certain algebraic properties, then all rings definable in ${\mathcal M}$ have some corresponding properties. Then, we notice that the same is true for groups and rings $*$-closed in a fixed compact $G$-space $(X,G)$.

\begin{theorem}\label{general 1}
(i) If every group definable in ${\mathcal M}$ is solvable-by-finite, then every ring with identity [or of finite characteristic] definable in ${\mathcal M}$ is nilpotent-by-finite.\\
(ii)  If every nilpotent group definable in ${\mathcal M}$  is abelian-by-finite, then every ring  definable in ${\mathcal M}$ is null-by-finite.\\
(iii) If every nilpotent group definable in ${\mathcal M}$ is (finite central)-by-abelian-by-finite, then every ring definable in ${\mathcal M}$ is (finite null)-by-null-by-finite.
\end{theorem}
{\em Proof.} (i) Let $R$ be a ring definable in ${\mathcal M}$ which contains 1 or which is of finite characteristic. Suppose $R$ has a finite characteristic $c$. Put $R_1=R \times {\mathbb Z}_c$, and define $+$ and $\cdot$ on $R_1$ by 
$$(a,k)+(b,l)=(a+b, k +_c l)\;\, \mbox{and}\;\, (a,k) \cdot (b,l)= (ab +k \times a + l \times b, k \cdot_c l).$$
This turns $R_1$ into a ring with 1 which is definable in ${\mathcal M}$ so that $R$ is a finite index ideal in $R_1$.
Thus, it is enough to consider the case when $R$ contains 1.

Let $H=Gl_3(R)$. Then, $H$ is definable in ${\mathcal M}$, and so it is solvable-by-finite. Let $H_0$ be a solvable, finite index subgroup of $H$. 

For distinct $i,j \in \{1,2,3\}$, define $t_{ij}(\alpha)$ as the element of  $H$ with 1's on the diagonal, $\alpha$ on the $(i,j)$-th position, and $0$'s elsewhere. We have the following well-known formula for commutators:
\begin{equation*}\tag{*}
[t_{ik}(\alpha),t_{kj}(\beta)]=t_{ij}(\alpha \beta),
\end{equation*}
for pairwise distinct $i,j,k \in \{1,2,3\}$.

Define $$I=\{ \alpha \in R: t_{ij}(\alpha) \in H_0\; \mbox{for all distinct}\; i,j \in \{1,2,3\} \}.$$
Using $(*)$ and the fact that $H_0$ is normal in $H$, we get that $I$ is an ideal of $R$. We claim that $I$ is of finite index in $R$. If not, then by Ramsey's theorem, there are $\alpha_i$,  $i\in \omega$, such that for some distinct $i,j \in \{1,2,3\}$ and for every $n>m$, we have $t_{ij}(\alpha_n - \alpha_m) \notin H_0$. But, $t_{ij}(\alpha_n - \alpha_m)=t_{ij}(\alpha_n)t_{ij}(\alpha_m)^{-1}$, and so $H_0$ is of infinite index in $H$, a contradiction. 

By the solvability of $H_0$, there is $n$ such that $H_0^{(n)}=\{ e\}$. $(*)$ implies that for any $\alpha_1,\dots,\alpha_{2^n} \in I$ and distinct $i,j \in \{ 1,2,3\}$, we have $t_{ij}(\alpha_1\dots \alpha_{2^n}) \in H_0^{(n)}= \{ e\}$, and so $\alpha_1\dots \alpha_{2^n}=0$. 

Thus, we have proved that $I$ is a nilpotent ideal of finite index in $R$.\\
(ii) %As in (i), it is enough to consider the case when $R$ is a ring interpretable in ${\mathcal M}$ which has an identity.
Let $H=UT_3(R)$ be the group of upper triangular $3\times3$ matrices with $1$'s on the diagonal (even if $1 \notin R$, we can use an external 1) and elements from $R$ above the diagonal. It is standard that $H$ is nilpotent. We see that $H$ is definable in ${\mathcal M}$, and so it is abelian-by-finite. Let $H_0$ be a normal, abelian subgroup of $H$ of finite index.

Define 
$$A= \left\{ (\alpha,\beta) \in R \times R:
\left(
\begin{array}{ccc}
1 & \alpha & \gamma\\
0 & 1 & \beta\\
0 & 0 & 1
\end{array}
\right) \in H_0\;\, \mbox{for some}\;\, \gamma \in R
\right\}.
$$

It is easy to check that $A$ is a finite index subgroup of $R^+ \times R^+$. Let $A_1$ denote the projection of $A \cap (R^+ \times \{ 0 \})$ on the first coordinate, and $A_2$ --  the projection of $A \cap (\{ 0 \} \times R^+)$ on the second coordinate. Then, $A_1$ and $A_2$ are finite index subgroups of $R^+$, and $A_1 \times A_2 <A$. Put $A_0=A_1 \cap A_2$. We see that $A_0$ is a finite index subgroup of $R^+$, and $A_0 \times A_0 <A$. We will be done if we show that $A_0$ is a null subring of $R$.

For any $\alpha,\beta,\gamma, \alpha',\beta',\gamma' \in R$ we have

$$ \left(
\begin{array}{ccc}
1 & \alpha & \gamma\\
0 & 1 & \beta\\
0 & 0 & 1
\end{array}
\right) 
\left(
\begin{array}{ccc}
1 & \alpha' & \gamma'\\
0 & 1 & \beta'\\
0 & 0 & 1
\end{array}
\right)
=
\left(
\begin{array}{ccc}
1 & \alpha + \alpha' & \gamma + \gamma' + \alpha \beta'\\
0 & 1 & \beta + \beta'\\
0 & 0 & 1
\end{array}
\right)
$$

Consider any $\alpha,\beta' \in A_0$. %Put $\alpha'=0$ and $\beta=0$. 
By the definition of $A$ and the fact that $H_0$ is abelian, there are $\gamma,\gamma' \in R$ such that

 $$ \left(
\begin{array}{ccc}
1 & \alpha & \gamma\\
0 & 1 & 0\\
0 & 0 & 1
\end{array}
\right) 
\left(
\begin{array}{ccc}
1 & 0 & \gamma'\\
0 & 1 & \beta'\\
0 & 0 & 1
\end{array}
\right)
=
\left(
\begin{array}{ccc}
1 & 0 & \gamma'\\
0 & 1 & \beta'\\
0 & 0 & 1
\end{array}
\right) 
\left(
\begin{array}{ccc}
1 & \alpha & \gamma\\
0 & 1 & 0\\
0 & 0 & 1
\end{array}
\right).
$$
This implies that $\alpha \beta'=0$. Thus, $A_0$ is a null subring.\\
(iii) We argue in a similar way as in (ii). Let $H=UT_4(R)$ be the group of upper triangular $4\times4$ matrices with $1$'s on the diagonal (as before, if $1 \notin R$, we can use an external 1) and elements from $R$ above the diagonal. Then, $H$ is nilpotent. We see that $H$ is definable in ${\mathcal M}$, and so it has a finite index, normal subgroup $H_0$ such that $[H_0,H_0]$ is finite and contained in the center of $H_0$. 

Define $A$ as the collection of all  $(\alpha_{12},\alpha_{23}, \alpha_{34}) \in R \times R \times R$ such that 
$$
\left(
\begin{array}{cccc}
1 & \alpha_{12} & \alpha_{13} & \alpha_{14}\\
0 & 1 & \alpha_{23} & \alpha_{24}\\
0 & 0 & 1 & \alpha_{34}\\
0 & 0 & 0 & 1
\end{array}
\right) \in H_0$$
for some $\alpha_{13},\alpha_{14},\alpha_{24} \in R$.

%Define $B$ as the collection of all $(\alpha_{13},\alpha_{24}) \in R \times R$ such that
%$$
%\left(
%\begin{array}{cccc}
%1 & 0 & \alpha_{13} & \alpha_{14}\\
%0 & 1 & 0 & \alpha_{24}\\
%0 & 0 & 1 & 0\\
%0 & 0 & 0 & 1
%\end{array}
%\right) \in H_0$$
%for some $\alpha_{14} \in R$.

One can check that $A$ is a finite index subgroup of $R^+ \times R^+ \times R^+$. 
%and $B$ is a finite index subgroup of $R^+ \times R^+$.
Let $A_1$, $A_2$ and $A_3$ be the projections of $A \cap (R \times \{ 0 \} \times \{0 \})$, $A \cap (\{0 \} \times R \times \{ 0 \})$ and $A \cap (\{ 0 \} \times \{ 0 \} \times R)$ on the first, on the second and on the third coordinate, respectively.  
%$B_1$ and $B_2$ are defined analogously. 
Let $A_0=A_1 \cap A_2 \cap A_3$. Then, $A_0$ is a finite index subgroup of $R^+$. Let $D$ be the subgroup of $R^+$ generated by $A_0$ and $A_0 \cdot A_0$. We will be done if we show that $D$ is a subring  of $R^+$ which is (finite null)-by-null. For this, it is enough to prove that the subgroup of $R^+$ generated by $A_0 \cdot A_0$ is finite and that $A_0 \cdot A_0 \cdot A_0=\{ 0 \}$.

Consider any $\alpha, \beta \in A_0$. Then, we can find
$$
M:=\left(
\begin{array}{cccc}
1 & \alpha & \alpha_{13} & \alpha_{14}\\
0 & 1 & 0 & \alpha_{24}\\
0 & 0 & 1 & 0\\
0 & 0 & 0 & 1
\end{array}
\right) \in H_0\;\;
\mbox{and}\;\;
N:=\left(
\begin{array}{cccc}
1 & 0 & \beta_{13} & \beta_{14}\\
0 & 1 & \beta & \beta_{24}\\
0 & 0 & 1 & 0\\
0 & 0 & 0 & 1
\end{array}
\right) \in H_0.$$
After some computations, we get

$$[M,N]=\left(
\begin{array}{cccc}
1 & 0 & \alpha \beta & \alpha \beta_{24}\\
0 & 1 & 0 & 0\\
0 & 0 & 1 & 0\\
0 & 0 & 0 & 1
\end{array}
\right) \in [H_0,H_0].$$

Since $[H_0,H_0]$ is finite, we get only finitely many possibilities for $\alpha\beta$, and so $A_0 \cdot A_0$ is finite. 
In fact, as all products of matrices (and their inverses) obtained above are still in $[H_0,H_0]$, we get that the subgroup of $R^+$ generated by $A_0 \cdot A_0$ is finite.

Consider any $\gamma \in A_0$. 
We can find $x,y,z \in R$ such that
$$
P:=\left(
\begin{array}{cccc}
1 & 0 & x & z\\
0 & 1 & 0& y\\
0 & 0 & 1 & \gamma\\
0 & 0 & 0 & 1
\end{array}
\right) \in H_0$$
We know that $[M,N] \in Z(H_0)$, so $P[M,N]=[M,N]P$. Computing these matrices, we get

$$
\left(
\begin{array}{cccc}
1 & 0 & \alpha \beta +x & \alpha \beta_{24} + z\\
0 & 1 & 0& y\\
0 & 0 & 1 & \gamma\\
0 & 0 & 0 & 1
\end{array}
\right) = 
\left(
\begin{array}{cccc}
1 & 0 & x +\alpha\beta& z + \alpha \beta \gamma + \alpha \beta_{24}\\
0 & 1 & 0& y\\
0 & 0 & 1 & \gamma\\
0 & 0 & 0 & 1
\end{array}
\right).$$
This implies that $\alpha \beta \gamma =0$, and so $A_0 \cdot A_0 \cdot A_0 = \{ 0 \}$.\hfill $\blacksquare$   

\begin{question}
Suppose that each solvable group definable in ${\mathcal M}$ is nilpotent-by-finite. Does it imply that every ring (with identity) definable in ${\mathcal M}$ is nilpotent-by-finite?
\end{question}

By the same arguments we get the following result.

\begin{theorem}\label{general 2}
Let $(X,G)$ be a compact $G$-space.\\  
(i) If every group $*$-closed in $X$ [or in $X^{teq}$] is solvable-by-finite, then every ring with identity [or of finite characteristic] $*$-closed in $X$ [$X^{teq}$] is nilpotent-by-finite.\\
(ii)  If every nilpotent group $*$-closed in $X$ [or in $X^{teq}$] is abelian-by-finite, then every ring $*$-closed in $X$ [$X^{teq}$] is null-by-finite.
\end{theorem}

Since each abelian, small compact $G$-group has finite exponent (see the proof of \cite[Proposition 2.3]{KW}), each small, compact $G$-ring has finite characteristic. Thus, assuming smallness in the above theorem, the extra assumption about the ring in point (i) is automatically  satisfied. In particular, we get

\begin{corollary}\label{general 3}
Let $(X,G)$ be a small profinite structure.\\  
(i) If every group $*$-closed in $X$ [or in $X^{teq}$] is solvable-by-finite, then every ring $*$-closed in $X$ [$X^{teq}$] is nilpotent-by-finite.\\
(ii)  If every nilpotent group $*$-closed in $X$  [or in $X^{teq}$] is abelian-by-finite, then every ring $*$-closed in $X$ [$X^{teq}$] is null-by-finite.
\end{corollary}

The following remark is rather standard.

\begin{remark}\label{definability for subgroups}
Let $G$ be any group.\\
(i) If $G$ has an abelian [normal] subgroup $H$ of finite index, then it has a definable, abelian [normal] subgroup of finite index which contains $H$.\\
(ii) If $G$ has a nilpotent [normal] subgroup $H$ of finite index, then it has a definable, nilpotent [normal] subgroup of finite index which contains $H$.
\end{remark}
{\em Proof.}
(i) %We can assume $H$ is normal. 
$C(H)$ contains $H$, so it is the intersection of finitely many centralizers of elements of $H$. Thus, $Z(C(H))>H$ is a definable, abelian [normal] subgroup of finite index in $G$.\\
(ii) Repeat the proof of \cite[Remark 3.3(ii)]{Kr3} using the fact that if $[G:C(N)]<\omega$, then $C(N)$ is definable (this allows one to eliminate the application of icc). Notice that when $G$ is countable,  $\omega$-categorical and $H$ is normal, the Fitting subgroup of $G$ (i.e. the group generated by all normal nilpotent subgroups) does the job. \hfill $\blacksquare$\\

Assuming icc on centralizers in definable quotients of definable subgroups, we have a variant of the above remark for solvability \cite[Remark 3.3(i)]{Kr3}. However, it is not clear to me what happens without icc.

\begin{question}
 Is it true that if $G$ has a solvable subgroup of finite index, then it has a definable, solvable, normal subgroup of finite index?
\end{question}

The answer to this question is positive for countable,  $\omega$-categorical groups: the product of all solvable, normal subgroups of finite index does the job.

We have the following counterpart of Remark \ref{definability for subgroups} for rings.

\begin{remark}\label{definability for subrings}
Let $R$ be any ring.\\
(i) If $R$ has a null ideal $S$ of finite index, then it has a definable, null ideal of finite index which contains $S$.\\
(ii) If $R$ has a nilpotent ideal $S$ of finite index, then it has a definable, nilpotent ideal of finite index which contains $S$.
\end{remark}
{\em Proof.}
(i) %By Remark \ref{Lewin}, we can assume $S$ is an ideal. 
$Ann_R(S)$  is a (two-sided) ideal which contains $S$, and so it is the intersection of annihilators of only finitely many elements of $S$. Thus, it is definable. We conclude that $Ann_R(Ann_R(S))$ is a definable, null ideal containing $S$.\\
(ii) %By Remark \ref{Lewin}, we can assume that $S$ is an ideal. 
The proof will be by induction  on the nilpotency class of $S$. If $S$ is null, we are done by (i). For the induction step, let $n\geq 3$ be the nilpotency class of $S$, and let $I=Ann_R(S^{n-1})$ and $J=Ann_I(I)$. Then, $I$ and $J$ are ideals of $R$. Moreover, $I$ contains $S$, so $J \subseteq Ann_R(S)$.  Using additionally that $S$ is of finite index in $R$, we see that $I$ is definable (because it is the intersection of annihilators of finitely many elements of $S^{n-1}$), and so $J$ is definable as well. Now, $(S+J)/J \cong S/(S \cap J)$ is of nilpotency class smaller than $n$ (as $S^{n-1} \subseteq J$). Thus, by the induction hypothesis, there is a definable ideal $K$ of $R$ such that $K/J$ is nilpotent. We can assume that $K \subseteq I$. Since $J =Ann_I(I)$ is null, we see that $K$ is nilpotent. \hfill $\blacksquare$\\

The same proofs yield the following remark.

\begin{remark}\label{definability remarks}
1. Let $(H,G)$ be a  $G$-group.\\
(i) If $H$ has an abelian [normal] subgroup $A$ of finite index, then it has a $*$-closed, abelian [normal] subgroup of finite index which contains $A$.\\
(ii) If $H$ has a nilpotent [normal] subgroup $N$ of finite index, then it has a $*$-closed, nilpotent [normal] subgroup of finite index which contains $N$.\\
2. Let $(R,G)$ be a $G$-ring.\\
(i) If $R$ has a null ideal $S$ of finite index, then it has a $*$-closed, null ideal of finite index which contains $S$.\\
(ii) If $R$ has a nilpotent ideal $S$ of finite index, then it has a $*$-closed, nilpotent ideal of finite index which contains $S$.
\end{remark}

\section{Applications}

First, we will use Theorems \ref{general 1} and \ref{general 2} to get structural results about $\omega$-categorical rings
and on  small compact $G$-rings.

\begin{theorem}
An $\omega$-categorical ring $R$ satisfying NSOP (the negation of the strict order property) is nilpotent-by-finite.
\end{theorem}
{\em Proof.} By Fact \ref{Mac NSOP}, we know that each group definable in $R$ is nilpotent-by-finite, and thus, by Theorem \ref{general 1}(i), we get that $R$ is nilpotent-by-finite. \hfill $\blacksquare$\\

The above result is a strengthening of \cite[Corollary to Theorem 3.6]{BR} which says that $\omega$-categorical stable rings are nilpotent-by-finite. The proof in \cite{BR} is completely different from ours, and it uses the full NOP (the negation of the order property).  

It was proved in \cite{KW} that each $\omega$-categorical, supesimple ring is (finite-null)-by-null-by-finite. Using our Theorem \ref{general 1}(iii), this result follows immediately from the fact that $\omega$-categorical, supersimple groups are (finite central)-by-abelian-by-finite \cite{EW}.

Now, we will use Theorem \ref{general 2} to get new information about small compact $G$-rings.

\begin{theorem}\label{nm-stable rings}
(i) A small, $nm$-stable compact $G$-ring is nilpotent-by-finite.\\
(ii) A small, compact $G$-ring of finite ${\cal NM}$-rank is null-by-finite.
\end{theorem}
{\em Proof.} 
(i) Use Fact \ref{nm-stable are solvable} and Theorem \ref{general 2}(i).\\
(ii) Use Fact \ref{NM finite implies abelian} and Theorem \ref{general 2}(ii). \hfill $\blacksquare$\\

\cite[Theorem 3.4]{KW} tells us that small, $nm$-stable profinite rings are null-by-finite. It is worth mentioning that this result is an immediate consequence of our Corollary \ref{general 3}(ii) and the fact that small, $nm$-stable profinite groups are abelain-by-finite \cite{Wa}. Similarly,  using \cite[Proposition 4.4]{Ne} which says that small, $m$-normal profinite groups are abelian-by-finite (see \cite{Ne1,Ne,Kr2} for the definition of $m$-normality), we get

\begin{corollary}\label{m-normal rings}
Each small, $m$-normal profinite ring is null-by-finite.
\end{corollary}
  
Theorem \ref{nm-stable rings}(ii) can be also obtained by application of \cite[Proposition 2.8]{KW2}. In order to see this, first we need to prove a lemma, which in fact will allow us to strengthen slightly Theorem \ref{nm-stable rings}(ii).

\begin{lemma}\label{clopen annihilator}
Let $(R,G)$ be a small, compact $G$-ring. Assume there is an $nm$-generic $g$ such that the left annihilator of $g$ is open in $R$. Then $R$ is null-by-finite.
\end{lemma}
{\em Proof.} There is an open (two-sided) ideal $I$ of $R$ contained in $Ann_R^{left}(g)$. Let $H=Stab(I)<G$. Then, by \cite[Remark 5.10]{Kr}, $H$ is a clopen subgroup of $G$ and it has countable index in $G$. Let $\{ g_i : i \in \omega\}$ be a set of representatives of left cosets of $H$ in $G$. Then,
$$\bigcup_{i \in \omega}g_iHg =o(g) \subseteq_{nm} R.$$ 
Thus, there is $i\in \omega$ for which $g_iHg \subseteq_{nm} R$, and so $Hg \subseteq_{nm} R$. But, $Hg \subseteq Ann_R^{right}(I) \lhd R$. Hence, $Ann_R^{right}(I)$ is clopen in $R$. This implies that $I \cap Ann_R^{right}(I)$ is a clopen, null ideal of $R$. \hfill $\blacksquare$   

\begin{theorem}
If $(R,G)$ is a small compact $G$-ring, and either ${\cal NM}(R)<\omega$ or ${\cal NM}(R) = \omega^{\alpha}$ for some ordinal $\alpha$, then $R$ is null-by-finite.
\end{theorem}
{\em Proof.} 
In the case ${\cal NM}(R)< \omega$, the idea of the proof is the same as in \cite[Theorem 3.4]{KW}. Take an $nm$-generic orbit $o(a/\emptyset)$. For $g \in o$, define $R_{g} = \{ (r,rg): r \in R\}$ -- a $g$-closed subgroup of $(R^+)^2$. By \cite[Proposition 2.8]{KW2}, there is $b \nmind a$ and a $b$-closed subgroup $K$ of $(R^+)^2$ commensurable with $R_a$. Take $a' \in o(a/b)$ such that $a' \nmind_b a$. Then, $R_a \sim K \sim R_{a'}$, so $R_a \sim R_{a'}$ (where $\sim$ denotes commensurability). Let $S$ be the projection of $R_a \cap R_{a'}$ on the first coordinate. Then, $S$ is a finite index subgroup of $R^+$ which is contained in $Ann_R^{left}(a-a')$. So, $Ann_R^{left}(a-a')$ is clopen in $R$. On the other hand, we see that $a' \nmind a$, and so $a-a'$ is $nm$-generic. We finish using Lemma \ref{clopen annihilator}. 

Now, consider the case ${\cal NM}(R)= \omega^{\alpha}$. In this case, we will often use Lascar inequalities for groups (see \cite[Fact 1.1]{KW2}). To be precise, we will be working in a slightly more general context in which $R$ is a $*$-closed ring in $X^{teq}$, where $(X,G)$ is a small compact $G$-space. By Theorem \ref{nm-stable rings}(i), $R$ is nilpotent-by-finite. Hence, by Remark \ref{definability remarks}, we can assume that $R$ is nilpotent. We argue by induction on the nilpotency class $n$ of $R$. If ${\cal NM}(Ann_R(R)) = \omega^{\alpha}$, then $Ann_R(R)$ has finite index in $R$, and we are done. If ${\cal NM}(Ann_R(R)) < \omega^{\alpha}$, then ${\cal NM}(R/Ann_R(R)) = \omega^{\alpha}$, so, by the induction hypothesis, we get that $R/Ann_R(R)$ is null-by-finite. By Remark \ref{definability remarks}, we can assume that it is null. Take an $nm$-generic $g \in R$, and consider a homomorphism $f_g: R^+ \to R^+$ given by $f_g(r)=rg$. Since $R/Ann_R(R)$ is null, $Im(f_g)\subseteq Ann_R(R)$. Thus, ${\cal NM}(ker(f_g)) =\omega^{\alpha}$, and so $Ann_R^{left}(g)$ is clopen in $R$. We finish using Lemma \ref{clopen annihilator}. \hfill $\blacksquare$ 

\begin{conjecture}
A small, $nm$-stable compact $G$-ring is null-by-finite. 
\end{conjecture}

A similar example to Example A of \cite{Kr} shows that the above conjecture strongly fails without the assumption of $nm$-stability.\\[3mm]
{\bf Example } Let $S_{\infty}$ be the group of all permutations of $\omega$. It acts on the ring $R:={\mathbb Z}_p^\omega$ permuting coordinates. Arguing as in \cite[Example A]{Kr}, we get that $(R,S_\infty)$ is a small, compact $S_\infty$-ring which is not nilpotent-by-finite.\\

Let us notice one more corollary of Lemma \ref{clopen annihilator}.

\begin{corollary}
Let $(R,G)$ be a small compact $G$-ring. If $R$ is countable-by-null-by-countable, then it is null-by-finite.
\end{corollary}
{\em Proof.} There is a countable index subring $R_1$ of $R$ and a countable ideal $I_1$ of $R_1$ such that $R_1/I_1$ is null. So, by the Baire category theorem, $R_1 \subseteq_{nm} R$. Hence, there is $g \in R_1$  which is $nm$-generic in $R$. We see that $Rg$ is countable, which implies that $Ann_R^{left}(g)$ is clopen in $R$. We finish using Lemma \ref{clopen annihilator}. \hfill $\blacksquare$\\  

The above results yield a quite good understanding of small, $nm$-stable compact $G$-rings.
As in the context of groups, a description of the structure of small [$nm$-stable] Polish $G$-rings (i.e. $G$-rings $(R,G)$, where $R$ is Polish) as well as searching examples of such rings are goals for the future.

In the next part of this section, using Theorem \ref{general 3}, we will obtain surprising relationships between some conjectures on small profinite groups. 

Recall the main conjecture on small profinite groups proposed by Newelski.

\begin{conjecture}\label{main profinite}
Each small profinite group is abelian-by-finite.
\end{conjecture}

The following three intermediate conjectures are still open.

\begin{conjecture}
For each small profinite group $(H,G)$:\\
(A) $H$ is solvable-by-finite,\\
(B) if $H$ is solvable-by-finite, then $H$ is nilpotent-by-finite,\\
(C) if $H$ is nilpotent-by-finite, then $H$ is abelian-by-finite.
\end{conjecture}

For small profinite rings, we have the following conjectures.

\begin{conjecture}
Each small profinite ring is null-by-finite.
\end{conjecture}

\begin{conjecture}
For each small profinite ring $(R,G)$:\\
(A') $R$ is nilpotent-by-finite,\\
(B') if $R$ is nilpotent-by-finite, then $R$ is null-by-finite.
\end{conjecture}

In fact, the proof of \cite[Remark 3.2]{KW} shows that (B') implies (A'). A step toward the proof of (A') was done in \cite{KW}, namely each small profinite ring was shown to be (nil of finite nil exponent)-by-finite. For us, \cite[Theorem 3.5]{KW} is essential.

\begin{fact}\label{ring to group}
Conjecture (A') (restricted to commutative rings) implies Conjecture (B). 
\end{fact}

Using this fact together with Corollary \ref{general 3}, we get the following

\begin{corollary}\label{last}
(i) (A) implies (A') implies (B).\\
(ii) (C) implies (B') implies (A') implies (B).
\end{corollary}

This corollary is surprising, because it implies that in order to show Conjecture \ref{main profinite}, it is enough to prove Conjectures (A) and (C) (although Conjecture (B) has not been proved). One could try to go further and to show that Conjecture \ref{main profinite} reduces to showing Conjecture (A). A possible way to do that could be through the following two conjectures.

\begin{conjecture}\label{last conjecture}
(i) (B) implies (B').\\
(ii) (B') implies (C).
\end{conjecture} 

%I do not see how to prove these conjectures.
% However, it is possible that Conjecture \ref{last conjecture}(i) can be proved by general arguments as in Section 1. 
An idea to prove Conjecture \ref{last conjecture}(ii) is by considering the so-called 'circle group of the ring' \cite{KP}. More precisely, in a small profinite group of nilpotency class 2 one should try to find a nilpotent ring structure which is $*$-closed and such that the original group is the circle group of this ring. Then, using (B'), it follows rather easily that the group we started from is abelian-by-finite. In fact, in the last implication, instead of $(B')$ it is enough to know that each small profinite ring is commutative-by-finite.   

Analyzing the proof of Fact \ref{ring to group}, one gets that all the results below Conjecture \ref{main profinite} are true in the class of small profinite groups and rings satisfying any extra assumption which is preserved under taking $*$-closed groups and rings in imaginary sorts. For example, $nm$-stability  and $m$-normality are preserved under taking $*$-closed groups and that is why \cite[Theorem 3.4]{KW} and our Corollary \ref{m-normal rings} follow from Corollary \ref{last} and the results saying that small, $nm$-stable [or $m$-normal] profinite groups are abelian-by-finite. 

%In particular, using \cite[Proposition 4.4]{Ne} which says that small, $m$-normal profinite groups are abelian-by-finite (see \cite{Ne1,Ne,Kr2} for the definition of $m$-normality), we get

%\begin{corollary}
%Each small, $m$-normal profinite ring is null-by-finite.
%\end{corollary}

While the general arguments from Section 1 transferring some properties of groups to the corresponding properties of rings are pretty easy, it is not clear what kind of converses are true. Fact \ref{ring to group} is an example of a result where some property of rings implies something about groups. But, the context in Fact \ref{ring to group}  is restricted to small profinite groups. We finish the paper with a similar result for $\omega$-categorical groups.

\begin{theorem}\label{last theorem}
Let ${\mathcal M}$ be an $\omega$-categorical structure. Suppose 
every group definable in ${\mathcal M}$ has a connected component and
every [commutative] ring interpretable in ${\mathcal M}$ is nilpotent-by-finite. Then every solvable group definable in ${\mathcal M}$ is nilpotent-by-finite.
\end{theorem}

Before the proof, recall that a group $G$ definable in a monster model $\C$ is said to have a connected component, $G^0$, if the intersection of all definable subgroups of finite index, denoted by $G^0$,  has bounded index (equivalently, it coincides with an intersection of only boundedly many definable subgroups of finite index, so it is type-definable). In the $\omega$-categorical context, this amounts to saying that $G^0$ is definable, or equivalently that there exists the smallest, definable subgroup of finite index. We take the last statement as the definition of the existence of connected components in non-saturated $\omega$-categorical groups.
 
Let us notice that the extra assumption about the existence of the connected components is satisfied for example under NIP. Indeed, if a monster model $\C \succ {\mathcal M}$ has NIP and $G$ is a group definable in $\C$, then $G^{00}$ (the smallest type-definable subgroup of bounded index) exists and is $\emptyset$-invariant, and so $\emptyset$-definable by $\omega$-categoricity. Thus, $G^{00}=G^0$. Moreover, in \cite{Kr4} we will see that $\omega$-categorical rings satisfying NIP are nilpotent-by-finite. So, Theorem \ref{last theorem} shows that solvable, $\omega$-categorical groups with NIP are nilpotent-by-finite.

Notice also that by Theorem \ref{last theorem}, all the relationships between conjectures on small profinite groups and rings formulated in Fact \ref{ring to group} and Corollary \ref{last} have their counterparts in the context of $\omega$-categorical structures satisfying any extra assumption which is preserved under interpretatability
and which implies (together with $\omega$-categoricity) that all definable groups have connected components.

Another remark is that from the proof below, it follows that in Theorem \ref{last theorem} it is enough to assume the existence of centralizer connected components (instead of connected components) for all interpretable groups. 
%Moreover, it is enough to assume that all commutative rings interpretable in ${\mathcal M}$ are nilpotent-by-finite.

Now, we turn to the proof of Theorem \ref{last theorem}.
In fact, we will modify the proof of \cite[Theorem 1.2]{Ma} (which says that solvable, $\omega$-categorical groups with NSOP are nilpotent-by-finite) replacing the part of the argument involving  NSOP by a new argument using virtual nilpotency of a certain interpretable ring. We could also give an argument based on the Baur-Cherlin-Macintyre proof that $\omega$-categorical, $\omega$-stable groups are abelian-by-finite (similarly to the proof of Fact \ref{ring to group}). We prefer to give a proof based on \cite{Ma}, because it leads to a question positive answer to which would allow us to drop the assumption about the existence of connected components in Theorem \ref{last theorem}.

First, recall Proposition 3.4 from \cite{Ma}, where we skip some unnecessary assumptions.

%\begin{fact}\label{Mac 1}
%If $G$ is a countable, $\omega$-categorical, metanilpotent  group, then all chief factors of $G$ are finite.
%\end{fact} 

\begin{fact}\label{Mac 2}
Let $G$ be a countable, $\omega$-categorical group with a normal subgroup $V$ which is a vector space over $F:=GF(p^a)$ ($p$ is a prime number). Let $H:=G/V$, and suppose $H$ has no elements of order $p$. Suppose $V$ is a sum of finite dimensional over $F$ (so finite) $FH$-modules. Then:\\
(i) $V$ is a direct sum of finite dimensional (over $F$), irreducible $FH$-modules,\\
(ii) $[H:C_H(V)]<\omega$. 
\end{fact}

Now, we will prove a suitable  variant of \cite[Corollary 3.5]{Ma}.

\begin{corollary}\label{Mac 3}
Let ${\mathcal M}$ be a countable, $\omega$-categorical structure and $G$ be a group interpretable in it.
Suppose $G$ has a normal subgroup $V$ interpretable in ${\mathcal M}$ which is a vector space over $F:=GF(p^a)$ ($p$ is a prime number). Assume that each ring interpretable in ${\mathcal M}$ is nilpotent-by-finite.  Let $H:=G/V$, and suppose $H$ is nilpotent and has no elements of order $p$.  Assume that each subgroup of $H$ interpretable in ${\mathcal M}$ has a connected component. Then $V$ is a direct sum of finite dimensional, irreducible $FH$-modules, and $[H:C_H(V)]$ is finite.
\end{corollary}
{\em Proof.}
We will argue by induction on the nilpotency class of $H$. First, we will reduce the induction step to the base step.
Suppose $H$ is of nilpotency class $n$, and the conclusion holds whenever the nilpotency class is smaller than $n$. Then, $Z_n(H)=H$. Put $H_1=Z_{n-1}(H)$. By the induction hypothesis, $[H_1:C_{H_1}(V)]$ is finite. Moreover, both $H_1$ and $C_H(V)$ are normal subgroups of $H$, so $C_{H_1}(V)$ is also normal in $H$. Since $[H,H] <H_1$, we conclude that all elements of $H/C_{H_1}(V)$ have centralizers of finite index. So, replacing $H$ by its connected component multiplied by $C_{H_1}(V)$, we can assume that $H/C_{H_1}(V)$ is abelian. Notice that it is enough to show the desired conclusion for $H/C_{H_1}(V)$ acting on $V$ by conjugation. Thus, everything has been reduced to the base step, i.e. to the case when $H$ is abelian. 

As in \cite{Ma}, we define $W$ as the sum of all finite dimensional $FH$-submodules of $V$. This is an interpretable in ${\mathcal M}$ subspace of $V$ invariant under $H$. By Fact \ref{Mac 2}, it is enough to show that $W=V$. Suppose for a contradiction that $W \subsetneq V$, and put $\overline{V}= V/W$.
Then, exactly as in \cite{Ma}, we get the following Claim.\\[3mm]
{\bf Claim} The $FH$-module $\overline{V}$ has no non-trivial, finite dimensional $FH$-submodules.\\[3mm]
The rest of the proof differs from \cite{Ma}. 
Choose a non-trivial $v \in \overline{V}$, and put $V_0=Lin_{F}(v^H)$ (the subgroup of $V$ generated by all elements $h^{-1}vh$, $h \in H$). By $\omega$-categoricity, $V_0$ is interpretable. Let $R$ be the ring of endomorphisms of $V_0$ generated by $H$. By the commutativity of $H$, we get that $R$ is interpretable in ${\mathcal M}$, and each $r \in R$ is determined by its value on $v$. So, $R$ has a nilpotent ideal $I$ of finite index $m$; say $r_1,\dots,r_m$ are representatives of all cosets of $I$ in $R$. By the claim, $V_0$ is infinite, and so $R$ and $I$ are infinite as well. Let $k \geq 2$ be the nilpotency class of $I$, i.e. the smallest number for which $I^k=\{0\}$. Choose any non-trivial $i \in I^{k-1}$. Then, $i(v) \ne 0$ and $i(v)^H \subseteq Ri(v) =\{r_1i(v),\dots,r_mi(v)\}$. Thus, $Lin_{F}(i(v)^H)$ is a non-trivial, finite dimensional (over $F$) $FH$-submodule of $\overline{V}$ the existence of which contradicts the claim. \hfill $\blacksquare$ \\

In order to prove Theorem \ref{last theorem}, we can assume that ${\mathcal M}$ is countable. Having Corollary \ref{Mac 3}, the proof of \cite[Theorem 1.2]{Ma} from page 490 of \cite{Ma} (which does not use NSOP anymore) goes through in our context and completes the proof of Theorem \ref{last theorem}.\\

The assumption about the existence of connected components in Theorem \ref{last theorem} was only used to reduce the induction step to the base step. We needed this reduction (to the case when $H$ is abelian) in order to make sure that the ring $R$ defined later in the proof  is interpretable. So, one can easily show that a positive answer to the following question would allow us to drop the extra assumption about connected components (assuming nilpotency of all rings interpretable in ${\mathcal M}$).

\begin{question}
We work in an $\omega$-categorical structure ${\mathcal M}$. Let $H$ be a definable group of nilpotency class at most 2 and with finite center. Suppose $H$ acts definably and by automorphisms on a definable,  abelian group $V$ [which is a vector space over $GF(p^a)$, and $H$ has no elements of order $p$]. Take $v \in V$, and put $V_0 = \langle Hv \rangle<V$. Let $R$ be the ring of endomorphisms of $V_0$ generated by $H$. Is it true that $R$ is interpretable in ${\mathcal M}$?
\end{question}

\noindent
Krzysztof Krupi\'nski\\
Instytut Matematyczny, Uniwersytet Wroc\l awski\\
pl. Grunwaldzki 2/4, 50-384 Wroc\l aw, Poland.\\
e-mail: kkrup@math.uni.wroc.pl\\[2mm]

\end{document}